\documentclass[12pt, a4paper, reqno, oneside]{amsart}
\usepackage{amssymb}  
\usepackage{amsmath} 
\usepackage{amsthm}
\usepackage{cite}

\newtheorem{theorem}{Theorem}
\newtheorem*{theorem*}{Theorem}

\newtheorem*{cor*}{Corollary}
\newtheorem{lemma}{Lemma}
\numberwithin{lemma}{section}
\newtheorem{prop}{Proposition}
\numberwithin{prop}{section}

\newtheorem*{remark*}{Remark}

\newtheorem*{definition*}{Definition}
\newtheorem{problem}{Problem}
\newtheorem*{problem*}{Problem}

\numberwithin{equation}{section}

\def\Op{\operatorname}
\def\mF{\mathbb{F}}
\def\mZ{\mathbb{Z}}
\def\cX{\mathcal{X}}
\def\cY{\mathcal{Y}}
\def\AGmL{\operatorname{A\Gamma L}}
\def\GmL{\operatorname{\Gamma L}}
\def\AGL{\operatorname{AGL}}
\def\PGmL{\operatorname{P\Gamma L}}
\def\PSL{\operatorname{PSL}}
\def\PGL{\operatorname{PGL}}
\def\SL{\operatorname{SL}}
\def\Soc{\operatorname{Soc}}
\def\Sym{\operatorname{Sym}}
\def\Alt{\operatorname{Alt}}
\def\Aut{\operatorname{Aut}}
\def\Iso{\operatorname{Iso}}
\def\ISO{\operatorname{ISO}}
\def\Orb{\operatorname{Orb}}
\def\Inv{\operatorname{Inv}}
\def\id{\operatorname{id}}
\def\poly{\operatorname{poly}}
\def\GL{\operatorname{GL}}
\def\TCL{\operatorname{TWOCLOSURE}}
\def\IMBED{\operatorname{IMBED}}

\def\BFC{\operatorname{BFC}}
\def\BFI{\operatorname{BFI}}
\def\AS{\operatorname{AS_0}}

\def\WL{\operatorname{WL}}
\def\32{$\frac{3}{2}$}
\def\12{$\frac{1}{2}$}
\def\2{{(2)}}

\begin{document}

\vspace{1cm}

\title[$2$-Closure of $\frac{3}{2}$-transitive group]{$\mathbf{2}$-Closure of $\mathbf{\frac{3}{2}}$-transitive group in polynomial time}

\author{Andrey V. Vasil$'$ev and Dmitry Churikov}

\thanks{The work was funded by RFBR according to the research project ¹ 18-01-00752}

\begin{abstract}
Let $G$ be a permutation group on a finite set $\Omega$. The $k$-closure $G^{(k)}$ of the group $G$ is the largest subgroup of $\Sym(\Omega)$ having the same orbits as $G$ on the $k$-th Cartesian power $\Omega^k$ of~$\Omega$. A group $G$ is called \32-transitive if its transitive and the orbits of a point stabilizer $G_\alpha$ on the set $\Omega\setminus\{\alpha\}$ are of the same size greater than one. We prove that the $2$-closure $G^{(2)}$ of a \32-transitive permutation group $G$ can be found in polynomial time in size of $\Omega$. In addition, if the group $G$ is not $2$-transitive, then for every positive integer $k$ its $k$-closure can be found within the same time. Applying the result, we prove the existence of a polynomial-time algorithm for solving the isomorphism problem for schurian \32-homogeneous coherent configurations, that is the configurations naturally associated with \32-transitive groups.

\end{abstract}

\maketitle

\section{Introduction}

Let $G$ be a permutation group on a finite set $\Omega$. The group $G$ is called \32-\emph{transitive} if it is transitive and orbits of the stabilizer $G_\alpha$ of a point $\alpha$ on the set $\Omega\setminus\{\alpha\}$ are of the same size greater than one. \32-Transitive groups naturally arise as normal subgroups of $2$-transitive groups. Examples of these groups are Frobenius groups and the automorphism groups of cyclotomic schemes over finite fields \cite{Delsarte,BCN} and near-fields \cite{BPR,VasChu}. The notion of a \32-transitive group was introduced by Wielandt \cite[\S~10]{WiFinPermG}, who laid the foundations of a theory of such groups by showing that every \32-transitive group is either primitive or Frobenius. Later Passman classified solvable \32-transitive groups \cite{Passman1, Passman2, Passman3}. Almost simple \32-transitive groups were described in~\cite{3/2AlmostSimple}, and the final step towards the classification of \32-transitive groups was done recently in~\cite{ArithmRes, 3/2}. The purpose of this paper is to show that this classification allows to prove that the $2$-closure problem for \32-transitive groups and the isomorphism problem for \32-homogeneous coherent configurations associated with these groups can be solved in polynomial time in degree of a group.

The key for our paper notion of $2$-closure (and wider, of $k$-closure) of a permutation group underlies Wielandt's method of invariant relations~\cite{WiInvRel}. The essence of the method consists of an application of the Galois correspondence between permutation groups on a finite set $\Omega$ and partitions of the $k$-th Cartesian power $\Omega^k$ of this set (see, e.\,g., \cite{FIK,EP_PermGroupApproach}). Every permutation group $G$ on  $\Omega$ is associated with the partition $\Orb_k(G)$ of $\Omega^k$ into the $k$-orbits of $G$, that is the orbits of the componentwise action of $G$ on $\Omega^k$. Conversely, every partition $P$ of the set $\Omega^k$ can be associated with the automorphism group $\Aut(P)$ of~$P$, which consists of all permutations on $\Omega$ that preserve classes of~$P$. Then inclusions
\begin{equation}\label{eq:Galois}
G\leq\Aut(\Orb_k(G))\text{ and }P\leq\Orb_k(\Aut(P)),
\end{equation}
express the correspondence.

According to~\cite{WiInvRel}, the group $G^{(k)}=\Aut(\Orb_k(G))$ from the first inclusion in~(\ref{eq:Galois}) is called the $k$-\emph{closure} of a group~$G$. Clearly, it is the largest subgroup of $\Sym(\Omega)$ having the same orbits as $G$ on~$\Omega^k$. The general problem we are interested in the present paper can be formulated as follows.

\medskip
\textbf{$k$-Closure problem.}
{\it Given a permutation group $G$ on a finite set $\Omega$ and a positive integer $k$, find the $k$-closure~$G^{(k)}$ of~$G$.}
\medskip

The case $k=1$ is not of particular interest, because $G^{(1)}$ is clearly the product of symmetric groups acting on orbits of $G$. The most important one is the $2$-closure problem which is equivalent to the problem of finding the automorphism group of some graph of initial group. Indeed, the pair $\Gamma=(\Omega,\Orb_2(G))$ can be considered as a complete colored graph on the vertex set $\Omega$, where every $2$-orbit corresponds to one of the colors. Then the group $G^{(2)}$ is the full automorphism group of this graph. The graph $\Gamma$ is a (colored) coherent configuration which is said to be schurian, because it associates with the group $G$ (see $\S$~\ref{sec:pcc} for details). In the sake of brevity we call such a configuration the \emph{scheme} of $G$ and denote it by~$\Inv(G)$. In this notation, the Galois correspondence from~(\ref{eq:Galois}) is expressed as
\begin{equation}\label{eq:Galois2}
G\leq\Aut(\Inv(G))\text{ and }\Gamma\leq\Inv(\Aut(\Gamma))
\end{equation}
(in the case $k=2$). The objects closed with respect to this correspondence are exactly $2$-\emph{closed} permutation groups, i.\,e., the groups with $G=G^{(2)}$, and \emph{schurian coherent configurations}, that is the coherent configurations on $\Omega$ satisfying the condition $\Gamma=\Gamma^{(2)}$, where $\Gamma^{(2)}=\Inv(\Aut(\Gamma))$.

In the context of computational complexity theory the above Galois correspondence leads to two natural problems: the finding of $2$-orbits of the automorphism group of a colored coherent configuration, i.e, given $\Gamma$, find $\Gamma^{(2)}$, and the $2$-closure problem. It is well-known that the former problem is polynomially equivalent to the general graph isomorphism problem and the later one can be polynomially reduced to it. In particular, this means that modulo the recent breakthrough result by Babai~\cite{Babai} both of this problems can be solved in time quasipolynomial in size of~$\Omega$. Yet, a class of permutation groups with known polynomial-time algorithms solving the $2$-closure problem is very restricted. Ponomarenko's article~\cite{Ponomarenko}, where the $2$-closure problem was considered from the computational complexity point of view, probably, at first time,  provided an algorithm polynomial in the degree of a group that solves this problem in the class of nilpotent groups. Later, the same result was obtained for the groups odd order~\cite{EvdokimovPonomarenko}. In the present paper we introduce a polynomial-time algorithm for finding the $2$-closure of an arbitrary \32-transitive permutation group.

\begin{theorem}\label{t:main}
The $2$-closure problem for a \32-transitive permutation group of degree $n$ can be solved in time polynomial in~$n$.
\end{theorem}

Let us discuss briefly this result here. First, note that we exploit standard polynomial-time algorithms from~\cite{Seress}. In particular, we assume that groups from input and output of our algorithm are given by generating sets of size polynomial in their degree. Second, the application of the classification of \32-transitive permutation groups  \cite[Corollary~3]{3/2} allows to prove that every such group which is not $2$-transitive either is $2$-closed, or is included in a group of order bounded by a polynomial of~$n$, so can be effectively constructed. Moreover, it turns out that obvious inclusions $G\leq G^{(k)}\leq G^{(2)}$ (with $k\geq2$) show that in this situation the general $k$-closure problem can be solved within the same time.

\begin{cor*}\label{cor}
Given a positive integer $k$, the $k$-closure problem for a \32-transitive but not $2$-transitive permutation group can be solved in time polynomial in its degree.
\end{cor*}

Finally, if the input group $G$ is $2$-transitive, then its $2$-closure, obviously, is the symmetric group, so is automatically done. The question how to find the $k$-closure of such a group for $k$ greater than $2$ remains open.

It is worth mentioning that in \cite{LPS,PS92} the problem of coincidence of the socle of a primitive group $G$ and the socle of its $k$-closure was studied. In particular, it was shown that $\Soc(G)=\Soc(G^{(k)})$ for all $k\geq6$. In the case where $G$ is a \32-transitive (not necessarily primitive) group, it is easy to deduce the same equality for all $k\geq2$ (see~$\S~2.2$).

Lets go back to a discussion on Theorem~\ref{t:main}. Recall that a coherent configuration is called \32-homogeneous if it is homogeneous and all its basis relations excluding the diagonal one have the same size. Thus the scheme of a \32-transitive group $G$ is exactly a schurian \32-homogeneous coherent configuration $\Inv(G)$. So Theorem~\ref{t:main} can be reformulated in combinatorial terms: the automorphism group of a schurian \32-homogeneous coherent configuration $\cX$ can be found in polynomial time (here we assume that the certificate of schurity of $\cX$ is given by a group $G$ with $\cX=\Inv(G)$). The natural continuation of this is the isomorphism problem for such configurations. More precisely, if we denote by $\psi$ a bijection between $\Orb_2(G)$ and $\Orb_2(G')$, and by $\Iso(\Inv(G),\Inv(G'),\psi)$ the set of all permutations $f$ in~$\Omega$ with $s^f=s^\psi$ for all $s\in\Orb_2(G)$, then the corresponding isomorphism problem can be formulated as follows.

\medskip
\textbf{Isomorphism problem for colored schurian coherent configurations.}
{\it Given permutation groups $G$ and $G'$ on a finite set $\Omega$ and a bijection $\psi$ between $\Orb_2(G)$ and $\Orb_2(G')$, find the set $\Iso(\Inv(G),\Inv(G'),\psi)$.}
\medskip

Clearly, if the set $\Iso(\Inv(G),\Inv(G'),\psi)$ is nonempty, then it is a coset of $G^{(2)}=\Aut(\Inv(G))=\Iso(\Inv(G),\Inv(G),\id)$ in~$\Sym(\Omega)$, here $\id$~is an identity map on $\Orb_2(G)$. Therefore, modulo the solved $2$-closure problem, in order to solve the isomorphism problem it suffices either to find a permutation from $\Sym(\Omega)$ which conjugates $G^\2$ and $(G')^{(2)}$, or to prove that there is no such permutation. It turns out that in the case of \32-transitive groups it can be done in polynomial time, so the following holds true.

\begin{theorem}\label{t:schur}
The isomorphism problem for colored schurian \32-homogeneous coherent configurations on $n$ points can be solved in time polynomial in~$n$.
\end{theorem}

In fact (see $\S$~\ref{sec:t2}), in the case when the corresponding coherent configurations are the schemes of primitive groups, we can find all isomorphisms between them (even those that do not preserve colors), i.\,e., we can find the set $$\Iso(\Inv(G),\Inv(G'))=\{g\in\Sym(\Omega)\mid s^g\in\Orb_2(G')\text{ for all }s\in\Orb_2(G)\}.$$

At the end of the introduction we wish to discuss the following open question. It is well known (see, e.\,g., \cite{EP_PermGroupApproach}) that there exist non-schurian coherent configurations (as well as there exist non-schurian  \32-homogeneous ones), so the schurity problem, including its complexity aspect, is one of the most substantial unsolved problem in theory of coherent configurations. In the context of this paper it is natural to draw attention to the particular case of this problem where a coherent configuration is supposed to be \32-homogeneous.

\begin{problem}\label{quest2} Is there an algorithm polynomial in $n$ which allows to determine whether or not a \32-homogeneous coherent configuration $\cX$ on $n$ points is schurian, and if so, to find a permutation group $G$ with $\cX=\Inv(G)$?
\end{problem}

In the next three sections we collect the necessary definitions, previously known results and auxiliary statements on permutation groups, coherent configurations and algorithmic tools. Sections 5 and 6 provide proofs of the main results.\medskip

Taking the opportunity the authors would like to express their deep gratitude to Ilia Ponomarenko for his valuable comments which allow to improve the presentation significantly. \medskip

\textbf{Notations.} Throughout the paper, $\Omega$~is a finite set of size~$n$.

For $\Delta\subseteq\Omega$ the diagonal of a Cartesian square $\Delta\times\Delta$ is denoted by $1_\Delta$, and for brevity, $1_\alpha=1_{\{\alpha\}}$.

For $s\subseteq\Omega\times\Omega$ put $s^*=\{(\beta,\alpha):\ (\alpha,\beta)\in s\}$.

If $\alpha\in\Omega$ and $s\subseteq\Omega\times\Omega$, then $\alpha s=\{\beta\in\Omega:\ (\alpha,\beta)\in r\}$.

If a group $G$ acts on $\Omega$, then for $g\in G$ and $\alpha\in\Omega$ denote by $\alpha^g$ the image of $\alpha$ under the action of~$g$. Set $s^g=\{(\alpha^g,\beta^g):\ (\alpha,\beta)\in s\}$ for $s\subseteq\Omega\times\Omega$.

The symmetric and alternating groups on the set $\Omega$ are denoted by $\Sym(\Omega)$ and $\Alt(\Omega)$ respectively.

Let $p$ be a prime, $d$ and $m$ positive integers, and $q=p^d$.

The cyclic group of order $m$, the elementary abelian $p$-group of order $p^d$ and the finite field of order~$q$ are denoted by $\mZ_m$, $\mZ_p^d$ and $\mF_q$ respectively.

We denote by $\GmL(m,q)$, $\GL(m,q)$, and $\SL(m,q)$ $m$-dimensional semilinear, linear and special linear group over the field $\mF_q$, and $\PGmL(m,q)$, $\PGL(m,q)$, and $\PSL(m,q)$ denote the projective versions of these groups.

We consider the affine group $\AGmL(m,q)$ ($\AGL(m,q)$ respectively) as a permutation group on $\Omega$ of the cardinality $q$ with a regular normal subgroup $V = \mZ_p^d$ where a stabilizer of a point acts on $V$ via conjugations as the group $\GmL(m,q)$ ($\GL(m,q)$ respectively).

\section{Preliminaries: permutation groups}\label{sec:pp}

In this section we compile necessary information on permutation groups (mostly about \32-transitive ones) and their $2$-closures. Here and further $G$ is a permutation group on a set~$\Omega$.

\subsection{\32-Transitive permutation groups}

We start with the classical Wielandt's result describing imprimitive \32-transitive groups.

\begin{lemma}\label{3/2Frobenius}\emph{\cite[Theorem~10.4]{WiFinPermG}}
A \32-transitive group is either primitive or Frobenius.
\end{lemma}

The next result reveals a general structure of primitive \32-transitive groups.

\begin{lemma}\label{3/2Primitive}\emph{\cite[Theorem~1.1]{3/2AlmostSimple}}
A primitive \32-transitive group is either affine of almost simple.
\end{lemma}

In the classification of solvable \32-transitive groups~\cite{Passman2}, there is an affine permutation group which we denote by $\AS(p^d)$ and call the Passman group. This group is of the form $\mZ_p^d \rtimes \Op{S}_0(p^{d/2})$, where $d$ is even, and $\Op{S}_0(p^{d/2})$ stands for the group of monomial matrices of dimension $2$ with determinant $\pm 1$:
$$
\Op{S}_0(p^{d/2}) = \left\langle\begin{pmatrix} \theta & 0 \\ 0 & \theta^{-1} \end{pmatrix}, \begin{pmatrix} -1 & 0 \\ 0 & 1 \end{pmatrix} , \begin{pmatrix} 0 & 1 \\ 1 & 0 \end{pmatrix}\right\rangle,
$$
and $\theta$ generates the multiplicative group $\mF_{p^{d/2}}^\times$\cite{Passman1}.

As already said, the classification of \32-transitive permutation groups was recently completed~\cite{3/2}. Below we summarize the main results from~\cite{3/2} for groups of sufficiently large degree.

\begin{lemma}\label{3/2}
Let $G$ be a \32-transitive permutation group of degree $n>169$. Then one of the following holds.
\begin{enumerate}

    \item $G$ is $2$-transitive.

    \item $G$ is Frobenius.

    \item $G$ is almost simple, $n = \frac{1}{2} q (q-1)$, $q = 2^f \geq 8$, $f$ is a prime; and either
    \begin{itemize}

        \item[(a)]$G = \PSL(2,q)$, and the size of a nontrivial orbit of a point stabilizer is equal to~$q+1$, or

        \item[(b)]$G = \PGmL(2,q)$, and the size of a nontrivial orbit of a point stabilizer is equal to~$f(q+1)$.

    \end{itemize}
    \item $G$ is affine, $n=p^d$, $p$ is a prime; and either
    \begin{itemize}

        \item[(a)]$G < \AGmL(1,p^d)$, or

        \item[(b)]$G = \AS(p^d)$ with $p$ odd and $d$ even.

    \end{itemize}
\end{enumerate}
\end{lemma}

\emph{Proof.} This follows from \cite[Corollaries 2 and~3]{3/2}. \qed\medskip

This classification yields the following assertions on the order and the minimal size of generating set of uniprimitive (i.\,e., primitive but not $2$-transitive) \32-transitive group.

\begin{lemma}\label{PolyOrder}
The order of an uniprimitive \32-transitive group $G$ of degree $n$ is polynomial in~$n$.
\end{lemma}

\emph{Proof.}
If $G$ is solvable, then its order is polynomially bounded by~\cite{Palfy}. Nonsolvable \32-transitive groups of degree greater than $169$ are Frobenius groups or projective linear groups of dimension $2$ due to Lemma~\ref{3/2}. In both cases the orders of the groups are polynomial in their degrees.\qed

\begin{lemma}\label{2gen}
A uniprimitive \32-transitive group $G$ of degree greater than $169$ is $2$-generated.
\end{lemma}

\emph{Proof.} Denote by $d(H)$ the minimal size of generating set of~$H$. By the main result of \cite{LM97}, if $N$ is a unique minimal normal subgroup of $G$, then $d(G)=\max\{2,d(G/N)\}$. Hence if $G$ is almost simple, then the required follows from Item 3 of Lemma~\ref{3/2}. Thus, applying Lemma~\ref{3/2} once again, we may assume that a minimal normal subgroup of $G$ is a regular elementary abelian $p$-group. Therefore, $d(G)=\max\{2,d(H)\}$, where $H$ is a point stabilizer. If $G$ is not Frobenius, then Item 4 of Lemma~\ref{3/2} implies that $H$ is contained in either $\GmL(1,p^d)$ or $\Op{S}_0(p^{d/2})$. Hence $H$ is a subgroup of a metacyclic group, i.\,e. it is generated by two elements. Finally, if $G$ is a Frobenius group, then the structure of its point stabilizer $H$ (so-called Frobenius complement) is well known (see, e.\,g., \cite[\S~18]{PassmanBook}). Applying \cite[Proposition~12.11]{PassmanBook}, one can easily deduce from~\cite[Theorems~18.2 and~18.6]{PassmanBook} that $d(H)\leq2$.  \qed

\subsection{$2$-Closures}

Recall that we denote the partition of the set $\Omega^k$ into $k$-orbits of $G\leq\Sym(\Omega)$ by $\Orb_k(G)$. Subgroups $G$ and $H$ of $\Sym(\Omega)$ are said to be $k$-equivalent if $\Orb_k(G)=\Orb_k(H)$. The largest $k$-equivalent to $G$ subgroup of $\Sym(\Omega)$ is called the $k$-closure of $G$ and denoted by $G^{(k)}$. In particular, $G\leq G^{(k)}$. A groups $G$ is called $k$-closed if $G = G^{(k)}$. Finally, a group $G$ is $k$-isolated if the $k$-equivalence of $G$ and $H$ implies $G=H$.

We note here the following basic fact, which, in particular, implies transitivity (primitivity) of the $k$-closure of a transitive (primitive) group.

\begin{lemma}\label{Inclusion}\emph{\cite[Theorem~5.7]{WiInvRel}}
If $H \leq G$, then $H^{(k)} \leq G^{(k)}$.
\end{lemma}

Further we focus on the case $k=2$. Obviously, a $2$-isolated group is $2$-closed. It turns out that the converse statement holds true for Frobenius groups.

\begin{lemma}\label{Frobenius}\emph{\cite[Lemma 4]{VasChu}}
Let the $2$-closure $G^\2$ of a group $G$ be a Frobenius group. Then $G=G^\2$.
\end{lemma}

Now we consider $2$-closures of \32-transitive groups. To begin with we note that the class of \32-transitive groups is closed with respect to taking $2$-closures

\begin{lemma}\label{3/2closure}
If $G$ is a \32-transitive group, then so is $G^\2$.
\end{lemma}

\emph{Proof.} As already mentioned, transitivity of $G$ yields transitivity of~$G^\2$. Due to a well-known correspondence between the $2$-orbits of a transitive group and the orbits of its point stabilizer, a coincidence of the $2$-orbits of groups $G$ and $G^\2$ implies a coincidence of the orbits of their point stabilizers. \qed\medskip

A solution of the $2$-closure problem for the imprimitive \32-transitive groups follows from the Lemma~\ref{3/2Frobenius} and the following statement.

\begin{lemma}\label{imp}\emph{\cite[Theorem~2.5.8]{FIK}}
An imprimitive Frobenius group is $2$-closed.
\end{lemma}

It follows from this lemma that imprimitive \32-transitive groups are, in fact, $2$-isolated. Another important subclass of \32-transitive groups also enjoys this property.

\begin{lemma}\label{3/2AlmostSimpleClos}
Let $G$ be a uniprimitive \32-transitive almost simple group of degree greater than~$169$. Then $G^\2 = G$.
\end{lemma}

\emph{Proof.} By above, $G^\2$ is uniprimitive and \32-transitive but not Frobenius. Due to Lemma~\ref{3/2}, the degree of $G$  equal to $n=\frac{1}{2} q (q-1)$, $q = 2^f \geq 8$, $f$  a prime, is not a prime power. Hence $G^\2$ cannot be affine. Finally, the groups $\PGmL(2,q)$ and $\PSL(2,q)$ cannot be $2$-equivalent because the sizes of nontrivial orbits of their point stabilizers are distinct (see Lemma~\ref{3/2} once again). \qed\medskip

The main assertion of this section summarizes the above information and gives a list of groups containing all \32-transitive groups of sufficiently large degree which are not $2$-closed.

\begin{prop}\label{p:not2closed}
Let $G$ be a \32-transitive group of degree greater than $169$. If $G\neq G^{(2)}$, then one of the following holds:
\begin{enumerate}
\item $G$ is $2$-transitive and $G^{(2)}=\Sym(\Omega);$
\item $G$ is primitive, $G<\AS(p^d)$, and $G^{(2)}=\AS(p^d);$
\item $G$ is primitive, $G<\AGmL(1,p^d)$, and $G^{(2)}<\AGmL(1,p^d).$
\end{enumerate}
\end{prop}

\emph{Proof.} If $G$ is an imprimitive group, then it is Frobenius due to Lemma~\ref{3/2Frobenius}, so by Lemma~\ref{imp} it is $2$-closed. The set $\Orb_2(G)$ of every $2$-transitive group consists of two elements, so $G^{(2)}=\Sym(\Omega)$. Thus we may assume that $G$ is uniprimitive. An application of Lemmas~\ref{Frobenius} and~\ref{3/2AlmostSimpleClos} shows that $G^\2$ satisfies the conditions of Item~4 of Lemma~\ref{3/2}, so the conclusion of the proposition follows. \qed\medskip

Now we return to the general case. The definition of the $k$-closure implies the following series of inclusions:
\begin{equation}\label{eq:chain}
G\leq\ldots\leq G^{(k+1)}\leq G^{(k)}\leq\ldots\leq G^{(2)}\leq G^{(1)},
\end{equation}
in particular, $G^{(k)}\leq G^\2$ for each $k\geq2$.

\begin{prop}\label{p:notKclosed} Let $k$ be a positive integer greater than~$1$, and $G$ a \32-transitive but not $2$-transitive group of degree greater than~$169$. If $G\neq G^{(k)}$, then $G$ is primitive and $G^{(k)}\leq G^{(2)}\leq H$, where $H\in\{\AS(p^d),\AGmL(1,p^d)\}$.
\end{prop}

\emph{Proof.} If $G\neq G^{(k)}$, then $G\neq G^{(2)}$ by the aforementioned remark. It follows from Proposition~\ref{p:not2closed} that $G$ is primitive and $G^{(k)}\leq G^\2\leq H$.\qed\medskip

At the end of this section we note (as already said in Introduction) that an easy consequence of the classification of \32-transitive groups and~\cite{BPR} is the following assertion on the coincidence of the socles of a \32-transitive group and its $k$-closure.

\begin{prop}\label{p:KcloserSoc} Let $k$ be a positive integer greater than $1$, and $G$ a \32-transitive but not $2$-transitive group. Then $\Soc(G)=\Soc(G^{(k)})$.
\end{prop}

\emph{Proof.} Due to Lemmas~\ref{3/2Frobenius} and~\ref{imp}, we may assume that $G$ is primitive. By Lemma~\ref{3/2Primitive}, it follows that $G$ is either almost simple or affine. Then we are done by \cite[Theorem~1.2]{3/2AlmostSimple} if $G$ is almost simple, and by~\cite[Theorem~3.2]{BPR} if it is affine. \qed

\section{Preliminaries: coherent configurations}\label{sec:pcc}

Here we collect the well-known facts on coherent configurations (see, e.\,g., \cite{EP_PermGroupApproach} and papers cited there).

\subsection{Main definitions}

Let $S$ be a partition of the set~$\Omega\times\Omega$ and $S^\cup$ the set of all unions of relations from~$S$. The pair $\cX=(\Omega,S)$ is called a \emph{coherent configuration} on~$\Omega$, if the following hold:
\begin{itemize}
\item[(C1)] $1_\Omega\in S^\cup$,
\item[(C2)] $s^*\in S$ for all $s\in S$,
\item[(C3)] for all $r,s,t\in S$, the number $c_{rs}^t=|\alpha r\cap\beta s^*|$
does not depend on the choice of $(\alpha,\beta)\in t$.
\end{itemize}
The elements of $\Omega$ and $S$, and the numbers $c_{rs}^t$ are called the {\em points} and {\em basis relations}, and the {\em intersection numbers} of the coherent configuration~$\cX$, respectively. The numbers $|\Omega|$ and $|S|$ are called the {\em degree} and the {\em rank} of~$\cX$.

A subset $\Delta$ of $\Omega$ is called a {\em fiber} of $\cX$, if $1_{\Delta}\in S$.  A coherent configuration~$\cX$  is called {\em homogeneous} if it has only one fiber, i.\,e., $1_\Omega\in S$. A homogeneous coherent configuration is \32-\emph{homogeneous} if for all $r,s\in S\setminus\{1_\Omega\}$ the equality $|r|=|s|$ holds.

A point $\alpha\in\Omega$ of a coherent configuration $\cX$ is said to be {\it regular}, if
$$
|\alpha r|\le 1\quad\text{for all}\ \,r\in S.
$$
If the set of regular points of a configuration $\cX$ is nonempty, then $\cX$ is called {\it 1-regular}; if it coincides with the set $\Omega$, then $\cX$ is called {\it semiregular}.

\subsection{Point extensions and the base size}

One can define a natural partial order on the set of coherent configurations over the same set by setting for $\cX=(\Omega,S)$ and
$\cX'=(\Omega,S')$:
$$
\cX\le\cX'\ \Leftrightarrow\ S^\cup\subseteq (S')^\cup.
$$
The minimal and maximal elements with respect to this order are the {\it trivial} and {\it complete} coherent configurations: the basis relations of the former are the diagonal $1_\Omega$ and its complement (for $n>1$), and all the basis relations of the latter are singletons.

Given two coherent configurations $\cX_1=(\Omega,S_1)$ and $\cX_2=(\Omega,S_2)$, there is the uniquely determined coherent configuration
$\cY=(\Omega,T)=\cX_1\cap\cX_2$ such that $T^\cup=(S_1)^\cup\cap(S_2)^\cup$. This allows us to define the {\em point extension} $\cX_{\alpha,\beta,\ldots}$ of a coherent configuration
$\cX=(\Omega,S)$ with respect to the points $\alpha,\beta,\ldots\,\in\Omega$ as follows:
$$
\cX_{\alpha,\beta,\ldots}=\bigcap_{\cY:\ S\subseteq T^\cup,1_\alpha,1_\beta,\ldots\in T^\cup}\cY,
$$
where $\cY=(\Omega,T)$. In other words, $\cX_{\alpha,\beta,\ldots}$ is the smallest coherent configuration on $\Omega$ which is larger or equal to $\cX$ and has the singletons $\{\alpha\},\{\beta\},\ldots$ as fibers.

A set $\Delta=\{\alpha,\beta,\ldots\}\subseteq\Omega$ is called a {\it base} of the coherent configuration $\cX$, if the extension $\cX_{(\Delta)}=\cX_{\alpha,\beta,\ldots}$ with respect to points $\alpha,\beta,\ldots$ from $\Delta$ is complete. The smallest cardinality of a base is called the {\it base number} of $\cX$ and denoted by $b(\cX)$ (sometimes, it is also, a little inaccurate, said to be the {\em base size} of~$\cX$). It is easy to see that the base size is bounded as follows: $0\le b(\cX)\le n-1$, and the equalities hold for the complete and trivial configurations respectively.

The next statement which can be considered as combinatorial analogous of Lemma~\ref{3/2Frobenius} plays a key role in considerations of the imprimitive case in the proof of Theorem~\ref{t:schur}.

\begin{lemma}\label{l:CombWielandtThm}
Let $\cX=\Inv(G)$ be a scheme of an imprimitive \32-transitive group~$G$. Then $b(\cX)=2$.
\end{lemma}

\emph{Proof.} This follows from~\cite[Theorem~5.11]{EP99_rus}. \qed

\subsection{Isomorphisms}

Coherent configurations $\cX=(\Omega,S)$ and $\cX'=(\Omega',S')$ are said to be {\it isomorphic}, if there exists a bijection $f:\Omega\to\Omega'$ such that the relation $s^f=\{(\alpha^f,\beta^f):\ (\alpha,\beta)\in s\}$ belongs to $S'$ for all $s\in S$. The bijection $f$ is called an {\it isomorphism} from $\cX$ onto $\cX'$; the set of all such isomorphisms is denoted by $\Iso(\cX,\cX')$. Obviously, the set $\Iso(\cX,\cX)$ is a permutation group on~$\Omega$.

Let now fix a bijection $\psi:S\to S'$. An isomorphism $f\in\Iso(\cX,\cX')$ with $s^f=s^\psi$ for all $s\in S$ is called an isomorphism of colored coherent configurations (w.r.t. $\psi$), and the set of all such isomorphisms is denoted by $\Iso(\cX,\cX',\psi)$. The subset $\Aut(\cX)=\Iso(\cX,\cX,\id_S)$ of elements of $\Iso(\cX,\cX)$, where $\id_S$ is an identity map on~$S$, is a normal subgroup of it and is called the {\it automorphism group} of (colored) coherent configuration~$\cX$.

A bijection $\varphi:S\to S',\ r\mapsto r'$, is called an {\it algebraic isomorphism} from~$\cX$ onto~$\cX'$ if
\begin{equation}\label{eq:separ}
c_{r^{}s^{}}^{t^{}}=c_{r's'}^{t'},\qquad r,s,t\in S.
\end{equation}
In this case coherent configurations $\cX$ and $\cX'$ are called {\it algebraically isomorphic}. Every isomorphism~$f$ from~$\cX$ onto~$\cX'$ naturally induces an algebraic isomorphism between them. The set of all isomorphisms inducing the algebraic isomorphism $\varphi$ is obviously equal to $\Iso(\cX,\cX',\varphi)$. If this set is nonempty for every algebraic isomorphism $\varphi:\cX\to\cX'$, then the coherent configuration $\cX$ is called {\it separable}. According to the formula~(\ref{eq:separ}), a configuration $\cX$ is separable if and only if it is determined by the collection of its intersection numbers up to isomorphism.

\begin{lemma}\label{l:semiregualr}
Every $1$-regular coherent configuration is separable. Moreover, if $\alpha$ and $\alpha'$~are regular points of coherent configurations $\cX$ and $\cX'$ lying in fibers $\Delta$ è $\Delta'$, respectively, and $\varphi$~is an algebraic isomorphism between these configurations such that $1_\Delta^\varphi=1_{\Delta'}$, then there exists a unique isomorphism $f$ from $\Iso(\cX,\cX',\varphi)$ mapping $\alpha$ to~$\alpha'$.
\end{lemma}

\emph{Proof.} The first statement is well known (see, e.\,g., \cite[Theorem~3.3]{EP_PermGroupApproach}), and the second one immediately follows from the definition of regular point. \qed

\section{Preliminaries: algorithms}\label{sec:pa}

In this section we describe algorithmic tools which allow to obtain the desired results.

\subsection{Standard algorithms}

As already said in Introduction, we assume that an input permutation group $G$ on $\Omega$ is given by its generating set of size polynomial in its degree~$n$. This means that in time $\poly(n)$ one can check whether the group $G$ is $2$-transitive, primitive (see, e.\,g., \cite{Seress}), as well as construct the scheme $\cX=\Inv(G)$ of~$G$.

We can check whether the cardinality $n$ of the set $\Omega$ is the power $p^d$ for a prime $p$, and if so, then construct a group $H\leq\Sym(\Omega)$ isomorphic to $\AGmL(1,p^d)$ or $\AS(p^d)$ with the socle $V$ acting regularly on~$\Omega$. In this case `to construct' means that we can list all elements of $H$ because the order of $H$ is clearly polynomial in~$n$.

If $H$ is a subset of $\Sym(\Omega)$ of cardinality polynomial in~$n$, then all automorphisms of a coherent configuration $\cX$ (all isomorphisms between configurations $\cX$ and $\cX'$) on~$\Omega$ lying in~$H$ can be found by brute force in time $\poly(n)$. Therefore, the sets $G^{(2)}\cap H$ and $\Iso(\Inv(G),\Inv(G'))\cap H$ can be found in polynomial time in~$n$. For convenience, we denote the results of the running of the last two algorithms by $\BFC(G;H)$ and $\BFI(\Inv(G),\Inv(G');H)$ respectively.

\subsection{Weisfeiler--Leman algorithm and its applications}

The classical Weisfeiler--Leman algorithm first appeared in~\cite{WL_NTI} and described in details in~\cite[Section~B]{W76} plays further a key role. The input of this algorithm is an arbitrary set $P$ of binary relations on a set  $\Omega$, and the output is the smallest coherent configuration
$$
\WL(P)=(\Omega,S)
$$
such that $P\subseteq S^\cup$. We call it the $\WL$-closure of~$P$. The running time of the algorithm is polynomial in the cardinalities of $P$ and~$\Omega$. An analysis of this algorithm in~\cite[Section~M]{W76} gives a proof of the following assertion (the exact statement is taken from~\cite[Theorem~2.4]{Ponomarenko2013}).

\begin{lemma}\label{l:wlai}
Let $P$ and $P'$ be two $m$-element sets of binary relations on an $n$-element set. Then given a bijection $\psi:P\to P'$,
one can check in time $mn^{O(1)}$ whether or not there exists an algebraic isomorphism $\varphi:\WL(P)\to\WL(P')$ such that $\varphi|_P=\psi$. Moreover, if $\varphi$ does exist, then it can be found within the same time.
\end{lemma}

Using this lemma one can easily show that there exists an efficient isomorphism-test algorithm between two algebraically isomorphic coherent configurations with the bounded base size. Namely, the following assertion holds.

\begin{lemma}\label{l:wlrb}\emph{\cite[Theorem~3.5]{Ponomarenko2013}}
Let $\cX$ and $\cX'$ be coherent configurations on $n$ points, and $\varphi:\cX\to\cX'$ an algebraic isomorphism between them. Then all the elements of the set $\Iso(\cX,\cX',\varphi)$ can be listed in time $(bn)^{O(b)}$, where $b=b(\cX)$ is the base number of~$\cX$.
\end{lemma}

Combining above two lemmas with Lemma~\ref{l:CombWielandtThm} we obtain a solution of the isomorphism problem for the colored schurian coherent configurations associated with imprimitive \32-transitive groups.

\begin{prop}\label{p:imprimitive}
Let $G$ and $G'$ be \32-transitive permutation groups on a set $\Omega$, $\psi$ a bijection between their sets of $2$-orbits. If the group $G$ is imprimitive, then the set $\Iso(\Inv(G),\Inv(G'),\psi)$ can be found in time $\poly(n)$.
\end{prop}

\medskip

The following well-known `group-theoretical' algorithm can be realized as `combinatorial' one with a help of the Weisfeiler--Leman algorithm.\medskip

Algorithm $\operatorname{\bf{IMBED}}$
\medskip

\textbf{Input:} a transitive group $G \leq \Sym(\Omega)$ with a $d$-element generating set $T$, a group $H \leq \Sym(\Omega)$ of order polynomial in~$n$, a $d$-element subset $T'$ of~$H$,
a bijection $\psi:T\to T'$, and a pair $(\omega,\omega')\in\Omega\times\Omega$.\medskip

\textbf{Output:} either a permutation $x$ on $\Omega$ such that $G^x \leq H$, $t^x=t^\psi$ for all $t\in T$, and $\omega^x=\omega'$, or the empty set if there is no such permutation.\medskip

\emph{Description of the algorithm}\medskip

\noindent{\bf Step 1:} Let $T=\{t_1,\ldots,t_d\}$, $T'=\{t_1',\ldots,t_d'\}$, and $t_i^\psi=t_i'$ for all $i=1,\ldots,d$. Define two sets $P=\{p_1,\ldots,p_d\}$ and $P'=\{p_1',\ldots,p_d'\}$ of binary relations on $\Omega$ as follows: for every $\alpha\in\Omega$ and $i=1,\ldots,d$, set $(\alpha,\alpha^{t_i})\in p_i$ and $(\alpha,\alpha^{t_i'})\in p_i'$. Fix a bijection $\psi:P\to P', p_i\mapsto p_i', i=1,\ldots,d$, induced by the bijection between sets $T$ and $T'$. \medskip

\noindent{\bf Step 2:} Applying the Weisfeiler--Leman algorithm, find the $\WL$-closures $\WL(P)$ and $\WL(P')$ of the sets $P$ and~$P'$. Note that by the definition of $P$ and $P'$ both coherent configurations $\WL(P)$ and $\WL(P')$ are semiregular. \medskip

\noindent{\bf Step 3:} Applying the algorithm from Lemma~\ref{l:wlai}, check whether or not there exists an algebraic isomorphism $\varphi$ between $\WL(P)$ and $\WL(P')$ such that $\varphi|_P=\psi$ and $1_\Delta^\varphi=1_{\Delta'}$, where $\Delta$ and $\Delta'$ are the fibers of coherent configurations $\WL(P)$ and $\WL(P')$ with $(\omega,\omega')\in\Delta\times\Delta'$.
\medskip

\noindent{\bf Step 4:} If the algebraic isomorphism $\varphi$ between $\WL(P)$ and $\WL(P')$ from Step 3 is found, then Lemma~\ref{l:semiregualr} and semiregularity of $\WL(P)$ and $\WL(P')$ imply that we can define a unique permutation $x\in\Iso(\WL(P),\WL(P'),\varphi)$ with $\omega'=\omega^x$. Put $$\IMBED((G,T),(H,T'),\psi,(\omega,\omega'))=x.$$

\noindent{\bf Step 5:} If there is no such algebraic isomorphism, then put $$\IMBED((G,T),(H,T'),\psi,(\omega,\omega'))=\varnothing.$$

\begin{prop}\label{p:imbed}
The algorithm $\IMBED$ is correct and runs in time~$\poly(dn)$.
\end{prop}

\emph{Proof.} If $x\in\Sym(\Omega)$ is such that $G^x \leq H$, $t^x=t^\psi$ for all $t\in T$, and $\omega'=\omega^x$, then it by construction lies in $\Iso(\WL(P),\WL(P'),\varphi)$, where $P$, $P'$, and an algebraic isomorphism $\varphi:\WL(P)\to\WL(P')$ are defined in the above algorithm. Vice versa, if $\Iso(\WL(P),\WL(P'),\varphi)\neq\varnothing$ and $1_\Delta^\varphi=1_{\Delta'}$ for the fibers $\Delta$ and $\Delta'$, then Lemma~\ref{l:semiregualr} and the semiregularity of $\WL(P)$ and $\WL(P')$ imply that there is $x\in\Iso(\WL(P),\WL(P'),\varphi)$ with $\omega'=\omega^x$. Then $$(\omega')^{t_i'}=(\omega^{t_i})^x=(\omega^x)^{x^{-1}t_i{x}}=(\omega')^{t_i^x}, i=1,\ldots,d.$$ Since $G$ is transitive, these equalities yield $t_i'=t_i^{x}$, $i=1,\ldots,d$. It remains to mention that the set $T$ generates the group $G$, and the set $T'$ generates a subgroup of~$H$ isomorphic to~$G$. Thus, the algorithm is correct. The bound on the running time of the algorithm follows from the bound on the running time of the Weisfeiler-Leman algorithm. \qed

\section{Proof of Theorem 1 and Corollary}\label{sec:t1}

The next algorithm solves the $2$-closure problem for the \32-transitive permutation groups.
\medskip

Algorithm $\operatorname{\bf{TWOCLOSURE}}$\medskip

\textbf{Input:} a \32-transitive permutation group $G$ on a set $\Omega$ of degree $n$.\medskip

\textbf{Output:} $G^\2=\TCL(G)$.\medskip

\noindent\emph{Description of the algorithm}

\medskip

\noindent{\bf Step 1:} If $n \leq 169$, put $\TCL(G)=\BFC(G;\Sym(\Omega))$. \medskip

\noindent{\bf Step 2:} If $G$ is $2$-transitive, put $\TCL(G)=\Sym(\Omega)$. \medskip

\noindent{\bf Step 3:} If $G$ is primitive and $n=p^d$ is a prime power, then fix a point $\omega\in\Omega$ and find all $2$-element subsets $T$ generating~$G$. \medskip

\noindent{\bf Step 4:} Construct $H\leq\Sym(\Omega)$ isomorphic to the Passman group $\AS(p^d)$. For every $2$-element subset $T$ generating $G$, every $2$-element subset $T'$ of $H$, every bijection $\psi:T\to T'$, and every point $\omega'\in\Omega$, apply the algorithm $$\IMBED((G,T),(H,T'),\psi,(\omega,\omega')).$$ If a permutation $x=\IMBED((G,T),(H,T'),\psi,(\omega,\omega'))$ is found, then put $$\TCL(G)=\BFC(G^x;H)^{x^{-1}}.$$

\noindent{\bf Step 5:} Construct $H\leq\Sym(\Omega)$ isomorphic to $\AGmL(1,p^d)$. For every $2$-element subset $T$ generating $G$, every $2$-element subset $T'$ of $H$, every bijection $\psi:T\to T'$, and every point $\omega'\in\Omega$, apply the algorithm $$\IMBED((G,T),(H,T'),\psi,(\omega,\omega')).$$ For all permutations $x=\IMBED((G,T),(H,T'),\psi,(\omega,\omega'))$, if they exist, find $$K(x)=\BFC(G^x;H).$$ Choose the permutation such that the order of $K(x)$ is maximal. Denote this permutation by~$y$, and put $\TCL(G)=K(y)^{y^{-1}}$.  \medskip

\noindent{\bf Step 6:} Put $\TCL(G)=G$.\medskip

The validity of Theorem~\ref{t:main} now follows from the validity of the following assertion.

\begin{prop}\label{p:tcl}
The algorithm $\TCL$ is correct and runs in time~$\poly(n)$.
\end{prop}

\emph{Proof.} The correctness of the algorithm follows from Propositions~\ref{p:not2closed} and~\ref{p:imbed}. Indeed, after Step 1 we may assume that the degree of $G$ is greater than 169. Suppose that $G\neq G^{(2)}$.

By Proposition~\ref{p:not2closed}, after Step 2 we may consider $G$ as a uniprimitive group of prime power degree~$p^d$. Lemma~\ref{2gen} implies that the group $G$ is $2$-generated, and this guarantees that the algorithm runs correctly at Step 3. Since $G$ is not $2$-transitive, either Item 2 of Proposition~\ref{p:not2closed} holds, i.\,e., $G\leq G^{(2)}\simeq\AS(p^d)$, or Item 3 does, i.\,e., $G\leq G^{(2)}\simeq K\leq\AGmL(1,p^d)$. Note that if  affine subgroups of $\Sym(\Omega)$, that is subgroups with a normal regular elementary abelian subgroup, are isomorphic, then are conjugated in~$\Sym(\Omega)$. The last fact follows from the conjugacy of isomorphic regular subgroups in the symmetric group and the equivalence of the natural action of a point stabilizer on $\Omega$ to its action by conjugation on the normal regular subgroup (see, e.\,g., \cite[Proposition~4.2]{PassmanBook}). Thus, if Item 2 of  Proposition~\ref{p:not2closed} holds, then there exists a permutation $x\in\Sym(\Omega)$ such that $(G^\2)^x=H$, where $H$ is the subgroup constructed at Step~4, and this permutation will be found via the algorithm $\IMBED$, whose correctness are provided by~Proposition~\ref{p:imbed}.

If Item 3 of Proposition~\ref{p:not2closed} holds, then among the permutations $x$ from $\Sym(\Omega)$ with $G^x\leq H\simeq\AGmL(1,p^d)$, where $H$ is the subgroup isomorphic to $\AGmL(1,p^d)$ and constructed at Step~5 (these permutations are again found by the algorithm $\IMBED$), we can choose a permutation $y$ such that $K=K(y)=(G^y)^{(2)}=(G^{(2)})^y\leq H$. Then for all other permutations $x$, we have $K(x)=(G^x)^{(2)}\cap H\leq K$, so $K$ is of the maximal order among all $K(x)$ that we found at Step~5.

Finally, if there are no permutations with a given property, then none of the items of Proposition~\ref{p:not2closed} holds, so the group $G=G^\2$ will be found at Step 6 of the algorithm.

Let us estimate the running time of the algorithm. The number of operations on Step 1, obviously, does not exceed a constant. The verification of $2$-transitivity at Step~2 as well as of primitivity at Step 3 is polynomial in~$n$. Lemma~\ref{PolyOrder} guarantees that the order of $G$, and so the number of $2$-element generating subsets~$T$ of $G$ is polynomial in~$n$. Finally, the polynomial-time algorithm $\IMBED$ is applied at Steps 4 and~5, and the number of calls of this algorithm is polynomial in the orders of $G$ and~$H$, which are, in turn, polynomial in~$n$. \qed
\medskip

It is easy to see why Corollary of Theorem~\ref{t:main} also holds true. Indeed, in the above algorithm we can omit Step~2. The rest follows from Proposition~\ref{p:notKclosed} and Theorem~\ref{t:main}.

\section{Proof of Theorem 2}\label{sec:t2}

Any nontrivial system of blocks of $G$ on $\Omega$ determines the equivalence relation on $\Omega\times\Omega$ distinct from $1_\Omega$ and~$\Omega^2$. Hence groups $G$ and $G^\2$ are primitive or imprimitive simultaneously. Therefore, due to Proposition~\ref{p:imprimitive}, in order to prove Theorem~\ref{t:schur} it suffices to deal with the case where $G$ and $G'$ are primitive. As already said in Introduction, we can prove even more in the case of primitive groups. Namely, applying the following algorithm $\operatorname{\textbf{ISO}}$ we can find the set $\Iso(\cX,\cX')$ of all (not only preserving colors) isomorphisms between the schemes $\cX=(\Omega,S)=\Inv(G)$ and $\cX'=(\Omega',S')=\Inv(G')$.\medskip

Algorithm $\operatorname{\textbf{ISO}}$
\medskip

\textbf{Input:} \32-transitive primitive permutation groups $G$ and $G'$ on a set $\Omega$ of degree $n$.\medskip

\textbf{Output:} the set $\Iso(\Inv(G),\Inv(G'))$. \medskip

\noindent\emph{Description of the algorithm}

\medskip

\noindent{\bf Step 1:} If $n \leq 169$, put $\Iso(\Inv(G),\Inv(G'))=\BFI(\Inv(G),\Inv(G');\Sym(\Omega)).$\medskip

\noindent{\bf Step 2:} Find $G^\2=\TCL(G)$ and $(G')^\2=\TCL(G')$.\medskip

\noindent{\bf Step 3:} If $G^\2=(G')^\2=\Sym(\Omega)$, put $\Iso(\Inv(G),\Inv(G'))=\Sym(\Omega)$. \medskip

\noindent{\bf Step 4:} If $|G^\2|=|(G')^\2|$, then fix $\omega\in\Omega$ and put $$\Iso(\Inv(G),\Inv(G'))=\{x\in\Sym(\Omega)\mid x=\IMBED((G^\2,T),((G')^\2,T'),\tau,(\omega,\omega'))\},$$ where $T$ runs over the $2$-element sets generating $G^\2$, $T'$ runs over the $2$-element sets generating $(G')^\2$, $\tau$ runs over the bijections from $T$ to $T'$, and $\omega'$ runs over~$\Omega$. \medskip

\noindent{\bf Step 5:} $\Iso(\Inv(G),\Inv(G'))=\varnothing$. \medskip

It remains to prove the following

\begin{prop}\label{p:iso}
The algorithm $\ISO$ is correct and finds the set $\Iso(\Inv(G),\Inv(G'))$ in time~$\poly(n)$.
\end{prop}

\emph{Proof.} As in the previous case we may assume that $n>169$. The groups $G^\2$ and $(G')^\2$ correctly and in polynomial time found at Step 2 due to Proposition~\ref{p:tcl} are the full automorphism groups of the schemes $\cX=\Inv(G)$ and~$\cX'=\Inv(G')$. If $G$ is $2$-transitive, then the scheme $\cX$ has rank~$2$, so it is isomorphic to the scheme $\cX'$ of $G'$ if and only if $G'$ is $2$-transitive. In this case every permutation from $\Sym(\Omega)$ is the required isomorphism. Thus, further we assume that $G$ and $G'$ are uniprimitive, hence the groups $G^\2$ and $(G')^\2$ are $2$-generated and their orders are polynomial in~$n$. In particular, it is easy to check the equality $|G^\2|=|(G')^\2|$, which is obviously the necessary condition for the set $\Iso(\cX,\cX')$ to be nonempty.

Note that the group $\Iso(\cX,\cX)$ is the normalizer in $\Sym(\Omega)$ of the group $\Aut(\cX)=G^\2$, and the same is true for the scheme $\cX'$ of the group~$G'$. Hence, given $x\in\Sym(\Omega)$, the equation  $(G^\2)^x=(G')^\2$ implies the equation $\Iso(\cX,\cX)^x=\Iso(\cX',\cX')$. As easily verified, the set of all permutations $x$ conjugating $\Iso(\cX,\cX)$ to $\Iso(\cX',\cX')$ coincides with the required set $\Iso(\cX,\cX')$. On the other hand, a permutation $x\in\Iso(\cX,\cX')$ moving $\cX$ to $\cX'$ conjugates the automorphism groups of these schemes, i.\,e., it conjugates $G^\2$ to~$(G')^\2$ in~$\Sym(\Omega)$. Thus, the required set coincides with the set of all permutations conjugating $G^\2$ to $(G')^\2$ which stays in the right-hand side of the equality at Step~4. It is found correctly and in polynomial time by Proposition~\ref{p:imbed} and the fact that sizes of the sets of all $T$, $T'$, $\tau$, and $\omega'$ defined at Step 4 are bounded by~$\poly(n)$. \qed\medskip

To complete the proof of Theorem~\ref{t:schur} it remains to note that having all isomorphisms between coherent configurations is easy to check which of them preserve colors.

\makeatletter
\renewcommand*{\@biblabel}[1]{\hfill\bf#1.}
\makeatother

\vspace{5mm}

Andrey V. Vasil$'$ev\smallskip

Sobolev Institute of Mathematics, 4 Acad. Koptyug avenue,

Novosibirsk State University, 1, Pirogova St.,

Novosibirsk, 630090, Russia\smallskip

{\em e-mail}: \verb"vasand@math.nsc.ru"

\vspace{5mm}

Dmitry Churikov\smallskip

Sobolev Institute of Mathematics, 4 Acad. Koptyug avenue,

Novosibirsk State University, 1, Pirogova St.,

Novosibirsk, 630090, Russia\smallskip

{\em e-mail}: \verb"churikovdv@gmail.com".

\end{document}